\documentclass{amsart}
\usepackage{graphicx}
\theoremstyle{definition}

\theoremstyle{remark}

\numberwithin{equation}{section}

\begin{document}

\title{Series Evaluation of a Quartic Integral}
\maketitle

\vskip 10pt

\begin{center}
{\bf Moa Apagodu}\\
{Department of Mathematics, Virginia Commonwealth University,\\ Richmond, VA 23284, USA}

\vskip 10pt

\end{center}

\begin{abstract}

  {We present a new single sum series evaluation of Moll's quartic integral and present two new generalizations. }

\end{abstract}

\vskip 10pt

\noindent In a beautiful personal story [6] Victor Moll describes his encounter with certain quartic integral and derives its evaluation and goes on to study analytic and number theoretic properties (\emph{log-concavity}, \emph{p-adic valuations}, \emph{location of the zeros}, etc.) of a polynomial associated with the evaluation of the integral [1,2,4,5,6,7]. In this article we use the Almkvist-Zeilberger algorithm ([3,8,9,10]) to derive a new series evaluation of this integral. In addition, we give two new generalizations of the identity. In [1], T. Amdeberhan and V. Moll presented a survey of old and new proofs of the evaluation and the formula: \\

\noindent \textbf{Theorem 1 [T. Amdeberhan and V. Moll, [1]]}:

$$
\int_{0}^{\infty}\frac{dx}{(x^4+2x^2a+1)^{m+1}}=\frac{\pi}{2}\frac{{2m\choose m}}{4^m(2(a+1))^{m+1/2}}{}_2F_{1}\left({{-m,m+1}\atop {-m+1/2}}\,;\,(a+1)/2\right)\,\,.
$$

\vskip .1in

\noindent where

$$
{}_2F_{1}\left({{a,b}\atop {c}}\,;\,x\right)=\sum_{k=0}^{\infty}\frac{(a)_k(b)_k}{(c)_k (1)_k}x^k\,\,
$$

\noindent and $(z)_k=z(z+1)(z+2)\ldots(z+k-1)$.\\

\noindent The polynomial associated with the evaluation of the integral that is the subject of study in [1,2,3,4] is

$$
P_m(a)=\frac{{2m\choose m}}{4^m}{}_2F_{1}\left({{-m,m+1}\atop {-m+1/2}}\,;\,(a+1)/2\right)\,.
$$

\noindent Next we state the main results of this article:\\

\noindent \textbf{Theorem 2}:

\begin{eqnarray*}
\int_{0}^{\infty}\frac{dx}{(x^4+2ax^2+1)^{m+1}}&=&\frac{1}{4}
\sum_{l=0}^{\infty}(-1)^l\frac{2^l(\frac{l}{2}-\frac{3}{4})!(m+\frac{l}{2}-\frac{1}{4})!}{l!m!}a^l\,\,.\\
\end{eqnarray*}

\noindent \textbf{Proof:}\\

\noindent  We use the Almkvist-Zeilberger algorithm ([3,8,9,10]), and the reader is assumed to be familiar with this method. In particular, we used Zeilberger's Maple package {\tt EKHAD8} (procedure AZc) that computes differential operators and certificates for single variable hyper-exponential functions accompanying [3], available from \\

{\tt http://www.math.rutgers.edu/\~{}zeilberg/tokhniot/EKHAD} .\\

\noindent We cleverly construct the (certificate) function

$$
R(x,a)=-\frac{x(4m+3+4ax^2m+2ax^2-x^4)}{(x^4+2ax^2+1)}
$$

\noindent with the motives

$$
-4m-3-4a(2m+3)D_a(F(x,a))-4(a^2-1)D_a^2F(x,a)=D_x(R(x,a)F(x,a))\,\,,
$$

\noindent where $F(x,a)$ is the integrand and $D_a$ is differentiation operator with respect to the variable $a$. If we integrate both sides with respect to $x$ on the limits of integration and observe that the right-hand side vanishes, we get a differential operator

$$
-4m-3-4a(2m+3)D_a-4(a^2-1)D_a^2\,\,,
$$

\noindent that annihilates the left side of the theorem. Using  the standard technique (or use Paul Zimmermann and Bruno Salvy's $gfun$ from Maple library if you wish) of translating a differential equation satisfied by a power series into a recurrence relation for its coefficients $a_l(m)$, we get

$$
(-4l^2+(-8m-8)l-4m-3)a_l(m)+(4l^2+12l+8)a_l(m)(l+2)=0\,\,,
$$

\noindent a homogeneous recurrence relation satisfied by the discrete coefficient function ${a_l(m)}$. Finally, the theorem follows by solving the recurrence relation with the initial conditions calculated directly from the integral: $a_0(m)=I(0,m)$ and $a_1(m)=I^\prime (0,m)$, where $I(a,m)$ is the the integral on the left. Q.E.D.\\

\bigskip

\noindent Comparing the right-hand side of our theorem with that of (\emph{theorem 1}), we get \\

$$
P_m(a)=
\frac{2^{m+3/2}(a+1)^{m+1/2}}{4\pi}\sum_{l=0}^{\infty}(-1)^l\frac{2^l(\frac{l}{2}-\frac{3}{4})!(m+\frac{l}{2}-\frac{1}{4})!}{l!m!}a^l\,\,\eqno{(Polypart)}\\
$$

\noindent Using Newton's Binomial theorem,

$$
(1+a)^{m+1/2}=\sum_{k=0}^{\infty}{m+1/2 \choose k}a^k\,\,,
$$

\noindent and multiplication of series, the coefficient of $a^n$, $d_n(m)$, in the polynomial $P_m(a)$ is \\

\begin{eqnarray*}
d_n(m)&=&\frac{2^{m+3/2}}{4\pi}\sum_{k+l=n} {m+\frac{1}{2} \choose k} (-1)^l\frac{2^l(\frac{l}{2}-\frac{3}{4})!(m+\frac{l}{2}-\frac{1}{4})!}{l!m!} \,\,\\
      &=&\frac{2^{m+3/2}}{4\pi}\sum_{l=0}^{n}{m+\frac{1}{2} \choose n-l} (-1)^l\frac{2^l(\frac{l}{2}-\frac{3}{4})!(m+\frac{l}{2}-\frac{1}{4})!}{l!m!}\,\,.
\end{eqnarray*}

\noindent Next, we give the first of two generalizations in which $2$ in the integral of \emph{theorem 2} is replaced by any integer $n$ for which the integral exists. \\

\noindent \textbf{Theorem 3}:

\begin{eqnarray*}
\int_{0}^{\infty}\frac{dx}{(x^{2n}+nax^{n}+1)^{m+1}}&=&
\frac{1}{2n}\sum_{l=0}^{\infty}(-1)^l\frac{n^l(\frac{l}{2}-\frac{2n-1}{2n})!(\frac{l}{2}+m-\frac{1}{2n})!}{l!m!}a^l\,\,.\\
\end{eqnarray*}

\noindent any integer $n$ for which the integral exists.\\

\noindent \textbf{Proof}:\\

\noindent First, we make the change of variables $z=x^n$ and the question reduces to evaluating

$$
\int_{0}^{\infty}\frac{dz}{n(z^2+2az+1)^{m+1}z^{1-1/n}} \,\,.
$$

\noindent Then, EKHAD gives a differential operator

$$
-2n-2nm+1 -(2m+3)n^2a D_a-n^2(a^2-1)D_a^2\,\,.
$$

\noindent with certificate function

$$
R(z,a)=-\frac{nz(2n+2nm-1+nzma+naz-az-z^2)}{z^2+az+1}\,\,.
$$

\noindent That is,

$$
(-2n-2nm+1 -(2m+3)n^2a D_a-n^2(a^2-1)D_a^2)F(x,a)=D_x(R(x,a)F(x,a))\,\,.
$$

\noindent where $F(x,a)$ is the integrand and $D_a$ is differentiation operator with respect to the variable $a$. Now integrate both sides and convert the resulting differential operator for the series into a recurrence relation for the coefficients and solve. Q.E.D.\\

\noindent The second generalization where $n$ is replaced by any parameter $\alpha$ for which the integral exists whose proof follows from \emph{theorem 3} by writing  $\alpha a $ as $n \left(\frac{a\alpha}{n}\right)$.\\

\noindent \textbf{Theorem 4}:

\begin{eqnarray*}
\int_{0}^{\infty}\frac{dx}{(x^{2n}+\alpha ax^{n}+1)^{m+1}}&=&
\frac{1}{2}\sum_{l=0}^{\infty}(-1)^ln^{l-1}(\frac{\alpha}{n})^{l}\frac{(\frac{l}{2}-\frac{2n-1}{2n})!(\frac{l}{2}+m-\frac{1}{2n})!}{l!m!}a^l\,\,.\\
\end{eqnarray*}

\noindent for any integer $n > 0$ and indeterminate $\alpha$ for which the integral exists.\\

\noindent \textbf{Problem}: Find analogous polynomial \emph{Polypart} as in theorem 1 associated with the evaluation of the generalization in \emph{theorem 3} if it exists, i.e.
find $closedForm(n,m,a)$ such that

\begin{eqnarray*}
P_m^n(a)&:=&
closedForm(n,m,a)\times\frac{1}{2n}
\sum_{l=0}^{\infty}(-1)^l\frac{n^l(\frac{l}{2}-\frac{2n-1}{2n})!(\frac{l}{2}+m-\frac{1}{2n})!}{l!m!}a^l\,\,,\\
\end{eqnarray*}

\noindent is a polynomial in $a$.\\

\noindent For the special case $n=2$ (\emph{Polypart}), $closedForm(2,m,a)=\frac{\pi}{2}\frac{1}{(2(a+1))^{m+1/2}}$. \\

\noindent {\bf References}\\

\bibliographystyle{amsplain}

\noindent [1]  T. Amdeberhan and V. Moll, A formula for a quartic integral: a survey of old proofs
and some new ones, Ramanujan J 18: 91–102 (2009).\\

\noindent [2]  T. Amdeberhan, D. Manna, and V. Moll, The 2-Adic Valuation of a Sequence Arising from a Rational Integral, Journal of Combinatorial Theory, Series A, \textbf{115}, 2008, 1474-1486.\\

\noindent [3]  M. Apagodu and D. Zeilberger, Multi-Variable Zeilberger and Almkvist-Zeilberger Algorithms
and the Sharpening of Wilf-Zeilberger Theory,  Adv. Appl. Math. 37(2006)(Special Regev issue), 139-152.\\

\noindent [4]  G. Boros, Victor H. Moll, and Sarah RileyAn,  An Elementary evaluation of a quartic integral,
Scientia, Series A: Math. Sciences \textbf{11}, 2005, 1-12.\\

\noindent [5]  M. Kauers and P. Paule, A Computer Proof of Moll's Log-Concavity Conjuctur, Proceedings of AMS, 135(12):3847--3856(2007).\\

\noindent [6]  V. Moll, The Evaluation of Integrals: a Personal Story, Notices Amer. Math. Soc. \textbf{49}, March 2002, 311-317.\\

\noindent [7] V. Moll and D. Manna, REMARKABLE SEQUENCE OF INTEGERS, to appear in Expositiones Mathematicae.\\

\noindent [8]  M. Petkov·sek, H.S. Wilf, D. Zeilberger, "A=B", A.K. Peters Ltd., 1996.\\

\noindent [9]  H. Wilf, D. Zeilberger, "Rational Functions Certify Combinatorial Identities", J. Amer. Math. Soc. 3 147-158 (1990). \\

\noindent [10]  D. Zeilberger, "The method of creative telescoping", J. Symbolic Computation, 11 195-204 (1991).

\end{document}